\documentclass{amsart}

\usepackage{amsmath}
\usepackage{amsthm}
\usepackage{amsfonts}
\usepackage{amssymb}
\usepackage{graphicx}

\newtheorem{thm}{Theorem}
\newtheorem{cor}{Corollary}
\newtheorem{lem}{Lemma}

\newtheorem{defn}{Definition}
\newtheorem{rem}{Remark}

\newcommand{\g}{{\gamma}}
\newcommand{\R}{{\mathbb{R}}}
\newcommand{\w}{{\omega}}

\begin{document}



\title{Intrinsically $n$-linked Complete Bipartite Graphs}

\author{Danielle O'Donnol}


\maketitle

\begin{abstract}
We prove that every embedding of $K_{2n+1,2n+1}$ into $\R^3$
contains a non-split link of $n$-components.  Further, given an
embedding of $K_{2n+1,2n+1}$ in $\R^3$, every edge of
$K_{2n+1,2n+1}$ is contained in a non-split $n$-component link in
$K_{2n+1,2n+1}$.
\end{abstract}



\section{Introduction}

A graph, $G$ is \emph{intrinsically linked} if every embedding of
$G$ into $\R^3$ contains a nontrivial link.  Conway and Gordon
\cite{CG} and Sachs \cite{Sa} first showed the existence of such
graphs by proving the complete graph on six vertices, $K_6$, is
intrinsically linked.  Sachs \cite{Sa} proved that the graphs in
the Petersen family are intrinsically linked and that no minor of
them is intrinsically linked. Then Robertson, Seymour, and Thomas
\cite{RST} proved that any intrinsically linked graph contains a
graph in the Petersen family as a minor. Together these results
fully characterize intrinsically linked graphs.

The idea of intrinsically linked graphs can be generalized to a
graph that intrinsically contains a link of more than two
components.  A link $L$ is \emph{split} if there is an embedding
of a 2-sphere $F$ in $\R^3\smallsetminus L$ such that  each
component of $\R^3\smallsetminus F$ contains at least one
component of $L$.  A graph $G$ is \emph{intrinsically $n$-linked}
if every embedding of $G$ into $\R^3$ contains a non-split
$n$-component link. Flapan, Naimi, and Pommersheim investigate
intrinsically triple linked graphs in \cite{FNP}.  They proved
that $K_{10}$ is the smallest complete graph to be intrinsically
triple linked. The more general question of what the smallest $m$
is such that $K_m$ is intrinsically n-linked is still open. Bowlin
and Foisy \cite{BF} also look at intrinsically triple linked
graphs. They exhibit two different subgraphs of $K_{10}$ that are
also intrinsically triple linked, proving that $K_{10}$ is not
minor minimal with respect to being intrinsically triple linked.
However, it is not known if either of these subgraphs is minor
minimal. Flapan, Foisy, Naimi, and Pommersheim address the
question of minor minimal intrinsically $n$-linked graphs in
\cite{FFNP}, where they construct families of minor minimal
intrinsically $n$-linked graphs.

In this paper we consider the question, for a given $n$, what is
the smallest $r$ and $s$ such that the complete bipartite graph
$K_{r,s}$ is intrinsically $n$-linked.  We prove that
$K_{2n+1,2n+1}$ is intrinsically $n$-linked.  A similar argument
shows that the complete tripartite graph $K_{2n, 2n, 1}$ is
intrinsically $n$-linked.

Though it is not known what the smallest $m$ is such that $K_m$ is
intrinsically $n$-linked, the number of vertices $m$ can be
bounded. First based on the number of disjoint simple closed
curves needed in $K_m$ we obtain a lower bound of $3n$.  This
lower bound is realized in the $n=2$ case, but no longer for the
$n=3$ case, where $m=10$ \cite{FNP}. It follows from \cite{FFNP}
that $K_{7n-6}$ is intrinsically $n$-linked.  Since $K_{2n,2n, 1}$
is a subgraph of $K_{2n+2n+1}$ Theorem \ref{thm2} gives an upper
bound for $m$ of $4n+1$.  This is a significant improvement over
the earlier bound of $7n-6$.

\section{Key Lemmas}

Let a simple closed curve containing exactly four vertices be
called a \emph{square}.  Let $\omega(J,K)=lk(J,K) \pmod{2}$. Hence
$\omega(J,K)$ is the number of times $J$ crosses over $K\pmod{2}$.
 When two simple closed curves, $J$ and $K$ are said to link, it
 should be understood that this means $\omega(J,K)\neq 0$.
 Let $K_{r,s}$ be embedded in $\mathbb{R}^3$, and let $\gamma$ be a
simple closed curve in $\mathbb{R}^3\smallsetminus K_{r,s}$. We
say $\gamma$ \emph{links} $K_{r,s}$ if there exists a square $J$
in $K_{r,s}$ such that $\omega(\gamma,J)\neq 0$.  The graph
$K_{3,3}$ has nine edges, and each edge is contained in four
squares.  Label by \{1,3,5\} and \{2,4,6\} the two sets of three
vertices of $K_{3,3}$. Each square can be uniquely denoted by a
four-number sequence, starting with an odd, and alternating odd
and even with the odds and evens each put in increasing order.

\begin{lem}\label{lemzero}
Let $K_{3,3}$ be embedded in $\mathbb{R}^3$, and let $\gamma$ be a
projection of a simple closed curve in $\mathbb{R}^3\smallsetminus
K_{3,3}$. Let $\gamma_o$ denote $\gamma$ with one of the
over-crossings in an edge $f$ changed to an under-crossing.
Suppose $\gamma_o$ links $K_{3,3}$ in zero squares.  Then $\gamma$
links $K_{3,3}$ in four squares, all containing the edge $f$.
\end{lem}

\begin{proof}
A square $J$ in $K_{3,3}$ links $\g$  if $\w(\g,J)\neq0$. Here
\[
\w(\g,J)=\w(\g_o,J)+
\begin{cases}
    1, & \text{if } f\in J \\
    0, & \text{if } f\not\in J \\
\end{cases}
\pmod{2}. \]  So \[ \w(\g,J)=0+
\begin{cases}
    1, & \text{if } f\in J \\
    0, & \text{if } f\not\in J \\
\end{cases}
\pmod{2} = \begin{cases}
    1, & \text{if } f\in J \\
    0, & \text{if } f\not\in J \\
\end{cases}
\pmod{2}. \]  The edge $f$, like any edge, is contained in four
squares in $K_{3,3}$.  Thus $\g$ links $K_{3,3}$ in four squares
all of which contain the edge $f$.
\end{proof}

\begin{lem}\label{lemfour}
Let $K_{3,3}$ be embedded in $\mathbb{R}^3$, and let $\gamma$ be a
projection of a simple closed curve in $\mathbb{R}^3\smallsetminus
K_{3,3}$. Let $\gamma_o$ denote $\gamma$ with one of the
over-crossings in an edge $f$ changed to an under-crossing.
Suppose $\gamma_o$ links $K_{3,3}$ in four squares, all of which
contain the edge e. Then:
\begin{itemize}
    \item if $f=e$, $\gamma$ links $K_{3,3}$ in zero squares
    \item if $f$ is adjacent to $e$, $\gamma$ links $K_{3,3}$ in four squares
    \item if $f$ is nonadjacent to $e$, $\gamma$ links $K_{3,3}$ in six
    squares.
\end{itemize}
\end{lem}

\begin{proof}
Again, if $J$ is a square in $K_{3,3}$, then
\[
\w(\g,J)=\w(\g_o,J)+ \begin{cases}
    1, & \text{if }f\in J \\
    0, & \text{if }f\not\in J \\
\end{cases}
\pmod{2}. \]  Here, \[ \w(\g_o,J)= \begin{cases}
    1, & \text{if }e\in J \\
    0, & \text{if }e\not\in J \\
\end{cases}
\pmod{2}.\]  So if $e=f$, \[\w(\g,J)=
\begin{cases}
    1+1, & \text{if } e=f\in J \\
    0+0, & \text{if }e=f\not\in J \\
\end{cases}
\pmod{2}=0\pmod{2}.\]  So $\g$ links $K_{3,3}$ in zero squares.

Next, if $f$ is adjacent to $e$, then it appears in two of the
four squares containing $e$.  Say, $J_1,J_2,J_3$ and $J_4$ are the
squares that contain $e$.  Let $J_3, J_4, J_5$ and $J_6$ be the
squares that contain $f$.  If $J$ is a square that contains
neither $e$ nor $f$, then $\w(\g,J)=\w(\g_o,J)=0$. So we need only
consider the above squares $J_i$. Now,
$$\w(\g,J_i)=1+0\pmod{2}=1 \pmod{2}, \text{  if } i=1,2,$$
$$\w(\g,J_i)=1+1\pmod{2}=0 \pmod{2}, \text{  if } i=3,4,$$ and
$$\w(\g,J_i)=0+1\pmod{2}=1 \pmod{2}, \text{  if } i=5,6.$$  So $\g$ links $K_{3,3}$ in
the four squares $J_1, J_2, J_5,$ and $J_6$.

Finally, if $f$ is nonadjacent to $e$ then it appears in exactly
one of the four squares containing $e$.  Let $J_1,J_2,J_3,$ and
$J_4$ be the squares that contain $e$, and let $J_4, J_5, J_6,$
and $J_7$ be the squares that contain $f$.  By an argument similar
to the above, we see that $\g$ links $K_{3,3}$ in the six squares
$J_1, J_2, J_3, J_5, J_6,$ and $J_7$.
\end{proof}

\begin{lem}\label{lemsix}
Let $K_{3,3}$ be embedded in $\mathbb{R}^3$, and let $\gamma$ be a
projection of a simple closed curve in $\mathbb{R}^3\smallsetminus
K_{3,3}$. Let $\gamma_o$ denote $\gamma$ with one of the
over-crossings in an edge $f$ changed to an under-crossing.
Suppose $\gamma_o$ links $K_{3,3}$ in six squares.  Then:
\begin{itemize}
    \item if $f$ appears in three of the squares that $\gamma_o$
links, $\gamma$ links $K_{3,3}$ in four squares
    \item if $f$ appears in two of
the squares that $\gamma_o$ links, $\gamma$ links $K_{3,3}$ in six
squares.
\end{itemize}
\end{lem}

\begin{proof}
Let $J_1, J_2, J_3, J_4, J_5,$ and $J_6$ be the squares that
$\g_o$ links in $K_{3,3}$.  Suppose $f$ appears in three of these
squares.  Let $f$ be contained in $J_4, J_5, J_6$, and $J_7$.
Similar to the proof of Lemma \ref{lemfour}, $\g$ links $K_{3,3}$
in the four squares $J_1, J_2, J_3,$ and $J_7$.

Next, suppose $f$ appears in two of these squares.  Let $f$ be
contained in $J_5, J_6, J_7$, and $J_8$.  Then $\g$ links
$K_{3,3}$ in the six squares $J_1, J_2, J_3, J_4, J_7,$ and $J_8$.
\end{proof}
We use the above lemmas to prove the following lemma.
\begin{lem}\label{lembig}
Let $K_{3,3}$ be embedded in $\mathbb{R}^3$, and let $\gamma$ be a
simple closed curve in $\mathbb{R}^3 \smallsetminus K_{3,3}$. Then
one of the following holds:
\begin{description}
    \item [(a)] $\gamma$ links $K_{3,3}$ in zero squares.
    \item [(b)] $\gamma$ links $K_{3,3}$ in four squares all of which
    contain a common edge.
    \item [(c)] $\gamma$ links $K_{3,3}$ in six squares and there exists three
    mutually nonadjacent edges of $K_{3,3}$ which appear in only two of the six
    squares. (All other edges appear in precisely three of the six squares.)
\end{description}
\end{lem}

\begin{proof}
Fix a projection of $\g$ and $K_{3,3}$.  The proof is by induction
on the number of over-crossings of $\g$ with $K_{3,3}$.

Suppose $n=1$.  So $\g$ has one over-crossing and it occurs on one
edge of $K_{3,3}$, say $e$. If we change the over-crossing to an
under-crossing, we obtain $\g_o$ which links zero squares in
$K_{3,3}$.  So by Lemma \ref{lemzero}, $\g$ links all of the
squares that contain $e$, of which there are four, and $\g$
satisfies \textbf{(b)}.

Suppose that $\g$ crosses over $K_{3,3}$ $n+1$ times.  Define
$\gamma_o$ to be the simple closed curve $\gamma$ with one of the
over-crossings in edge $f$ changed to an under-crossing.  By our
inductive hypothesis, $\g_o$ satisfies \textbf{(a)}, \textbf{(b)},
or \textbf{(c)}.  If $\g_o$ satisfies \textbf{(a)} then by Lemma
\ref{lemzero}, $\g$ satisfies \textbf{(b)}. If $\g_o$ satisfies
\textbf{(b)}, we may assume without loss of generality, $\g_o$
links the squares 1234, 1236, 1254, and 1256, all containing the
edge $\overline{12}$.  Hence, from Lemma \ref{lemfour} we know
that $\g$ links zero, four, or six squares in $K_{3,3}$.  If the
edge $f$ is adjacent to $\overline{12}$, without loss of
generality $f=\overline{23}$. Then, by the proof of Lemma
\ref{lemfour}, $\g$ links the squares 1254, 1256, 3254, and 3256.
So $\g$ links four squares that all contain the edge
$\overline{25}$ and thus it satisfies \textbf{(b)}. Next, if the
edge $f$ is nonadjacent to $\overline{12}$, without loss of
generality $f=\overline{34}$. By the proof of Lemma \ref{lemfour},
$\g$ links the squares 1254, 1436, 1236, 3254, 1256,and 3456.  In
these six squares, the mutually nonadjacent edges $\overline{14}$,
$\overline{32}$, and $\overline{56}$ each appears only twice
(while all other edges of $K_{3,3}$ appear thrice).  So the curve
$\g$ satisfies \textbf{(c)}.

If $\g_o$ satisfies \textbf{(c)}, without loss of generality
$\g_o$ links the six squares where the edges $\overline{12}$,
$\overline{34}$, and $\overline{56}$ appear only twice.  Then
$\g_o$ links 1236, 1254, 1436, 3254, 1456, and 3256. Every edge of
$K_{3,3}$ either appears in two or three of the squares that
$\g_o$ links. If $f$ appears in three, without loss of generality
$f=\overline{14}$.  If $f$ appears in two, without loss of
generality $f=\overline{12}$. Recall that, $\g$ links a square $J$
if and only if either $J$ contains $f$ and $J$ does not link
$\g_o$, or $J$ does not contain $f$ and $J$ links $\g_o$.  If
$f=\overline{14}$, then $\g$ links the four squares 1234, 1236,
3254, and 3256, which all contain the edge $\overline{23}$. Hence
$\g$ satisfies \textbf{(b)}.  If $f=\overline{12}$, then $\g$
links the six squares 1234, 1256, 1436, 3256, 1456, and 3254, in
which the three mutually nonadjacent edges $\overline{12}$,
$\overline{36}$, and $\overline{54}$ appear only twice.  Hence
$\g$ satisfies \textbf{(c)}.
\end{proof}

From here forward, when it is said that $\g$ links a graph
isomorphic to $K_{3,3}$ in four squares, it should be understood
that $\g$ satisfies \textbf{(b)} in Lemma \ref{lembig}, as this is
the only way $\g$ can link $K_{3,3}$ in four squares. Similarly,
if $\g$ is said to link six squares in $K_{3,3}$, that should be
taken to mean $\g$ satisfies \textbf{(c)}.

\section{Main results}

We will use the following definition in proving the main result.
\begin{defn}
Let $M$ be a subgraph of $K_{4,4}$ which is isomorphic to
$K_{3,3}$, and let $\alpha$ be a square in $M$.  The
\emph{$\alpha$-opposite subgraph of $M$}, say $N$, isomorphic to
$K_{3,3}$ is the subgraph defined by the four vertices in $\alpha$
and the two vertices in $K_{4,4}\smallsetminus M$.  Then $M$ and
$N$ are a \emph{pair of $\alpha$-opposite subgraphs} in $K_{4,4}$.
\end{defn}

\begin{thm}\label{Thm1}
Given $n>1$, every embedding of the complete bipartite graph
$K_{2n+1,2n+1}$ into $\mathbb{R}^3$ contains a non-split
$n$-component link.
\end{thm}

\begin{proof}
We shall prove by induction on $n$ that every embedding of
$K_{2n+1,2n+1}$ contains a non-split $n$-component link $L$ of
squares, with a component $J$ such that $L\smallsetminus J$ is a
non-split $(n-1)$-component link.  When $n=2$,
$K_{2n+1,2n+1}=K_{5,5}$. The graph $K_{4,4}\subset K_{5,5}$ is
known to be intrinsically 2-linked (by Sachs \cite{Sa}, $K_{4,4}$
with one edge removed is intrinsically linked). In this case,
since $K_{4,4}$ contains a non-split 2-component link, both
components must be squares and either component can be chosen to
be $J$.

Consider an embedding of $K_{2n+1,2n+1}$ into $\mathbb{R}^3$. Fix
a projection of the embedded $K_{2n+1,2n+1}$.  It has a subgraph,
say $H$, which is isomorphic to $K_{2(n-1)+1,2(n-1)+1}$. By the
inductive hypothesis, $H$ contains a non-split $(n-1)$-component
link $L$ of squares with a component $J$ such that
$L\smallsetminus J$ is a non-split $(n-2)$-component link. Let
$H_1$ be the subgraph of $H$ that is defined by the vertices of
$L$.  Next choose a square $\g$ in $L$ that links $J$.   Let the
subgraph isomorphic to $K_{5,5}$, defined by the vertices of
$(K_{2n+1,2n+1}\smallsetminus H_1)\cup J$ be called $G$. Since
$\g$ is in $H_1$ and is disjoint from $J$, it is disjoint from
$G$. The curve $\gamma$ links $J$ so, by Lemma \ref{lembig},
$\gamma$ links each of the $K_{3,3}$ subgraphs containing $J$ in
four or six squares.

We will consider two different cases which each break into two
subcases and show in each case that there are two disjoint squares
in $G$ that either both link $\gamma$ or they are linked together
and one of them links $\g$. This will finish the proof by finding
the desired non-split $n$-link in the given embedding of
$K_{2n+1,2n+1}$.
\medskip

\textbf{Case 1:} The curve $\g$ does not link any of the $K_{3,3}$
subgraphs of $G$ in six squares.  Thus every $K_{3,3}$ subgraph of
$G$ that links $\g$, links $\g$ in four squares which all share a
common edge.
\smallskip

Let $G_0$ be a subgraph of $G$ isomorphic to $K_{4,4}$ containing
$J$.  Label the two set of vertices of $G_0$ by \{1,3,5,7\} and
\{2,4,6,8\}, and the remaining vertices of $G$ label 9 and 0
appropriately.
\smallskip

\textbf{Case 1(a):}  For every square $\alpha$ in $G_0$ and edge
$e$ of $\alpha$, there is no pair of $\alpha$-opposite subgraphs
of $G_0$ such that $e$ is the common edge of the four squares
linking $\g$ in both subgraphs.
\smallskip

Suppose without loss of generality, that $\g$ links 123456 in the
four squares 1234, 1236, 1254, and 1256 (common edge
$\overline{12}$), and links the 1234-opposite subgraph 123478 in
four squares with a different common edge.  Since $\g$ links 1234
in 123478 the common edge must be an edge of 1234.  There are two
different ways this can happen, either the common edge is adjacent
to $\overline{12}$ or not.  First consider the linking where the
common edge ($\overline{34}$) is not adjacent to the first common
edge ($\overline{12}$). Then $\g$ links 123478 in the four squares
1234, 1438, 3274, and 3478.  So $\g$ links the two disjoint
squares 1256 and 3478 which are both in $G_0$.  Call 1256, $L_o$
and 3478, $L_1$ . Now $(L\smallsetminus J)\cup L_o\cup L_1$ is a
non-split $n$-component link, and $(L\smallsetminus J)\cup L_o$ is
a non-split $(n-1)$-component link.

Next suppose that $\g$ links 123478 in the four squares with
common edge adjacent to $\overline{12}$, say $\overline{14}$. So
$\g$ links 1234, 1274, 1438 and 1478 in 123478. We see as follows
that this forces $\g$ to link twelve squares in $G_0$. The simple
closed curve $\g$ links the three squares 1234, 1438 and 1254 in
123458. These squares have a single edge $\overline{14}$ in
common, so $\g$ must also link 1458. Similarly, $\g$ links 1234,
1236, and 1274 in 123476 and therefore the additional square 1276.
Next, $\g$ links the two squares 1236 and 1256 in 123658. These
squares have two edges in common $\overline{12}$ and
$\overline{16}$. However $\overline{12}$ cannot be the common edge
because the squares 1258 and 1238 are in the subgraph 123458, and
are not among the four squares $\g$ links in this subgraph. Thus,
$\overline{16}$ is the common edge for 123658.  So $\g$ also links
1638 and 1658. Finally, $\g$ links the three squares 1438, 1458
and 1638 in 143678, and therefore also links 1678.  So we have
found twelve squares of $G_0$, 1234, 1236, 1254, 1256, 1274, 1276,
1438, 1638, 1458, 1658, 1478, and 1678 that $\g$ links. By
inspection the $K_{3,3}$ subgraphs of $G_0$: 125478, 125678, and
145678 as well as the above mentioned $K_{3,3}$ subgraphs of
$G_0$: 123456, 123476, 123458, 123478, 123658, and 143678 each
contains exactly four of the above mentioned twelve squares (that
$\g$ links) with a common edge.  Together these nine $K_{3,3}$
subgraphs of $G_0$ contain all of the squares in $G_0$.  Thus
these are the only squares in $G_0$ that $\g$ links, because in
this way we see that it does not link any other square in $G_0$.
If none of these twelve squares in $G_0$ is contained in a link in
$G_0$, then there is not a pair of squares in $G_0$ that will form
a non-split $n$-component link together with $L$.

Let the subgraph of $G$ defined by the vertices 14365870 be $G_1$.
The simple closed curve $\g$ links 143658 in the four squares
1438, 1458, 1638, and 1658.  Notice the common edge is
$\overline{18}$.  The simple closed curve $\g$ also links the
1438-opposite subgraph 143870. There are three different ways this
can happen. If $\g$ links 143870 with common edge $\overline{18}$
then we are in \textbf{case 1(b)}.  If $\g$ links 143870 with the
common edge $\overline{34}$ then there are two disjoint squares
1658 and 4307 linking $\g$ as described two paragraphs above. If
$\g$ links 143870 in four squares with common edge
$\overline{14}$, then it links twelve squares in $G_1$, which we
can find as before: 1438, 1458, 1478, 1638, 1658, 1678, 1430,
1450, 1470, 1630, 1650, and 1670. Again there is a possibility
that the links in $G_1$ do not contain any of these squares.
Finally, consider the subgraph 12365870.  The simple closed curve
$\g$ links 123658 in the four squares 1236, 1256, 1638, and 1658
(common edge $\overline{16}$) and it links the 1236-opposite
subgraph 123670 in the four squares with the same common edge
$\overline{16}$ (i.e. 1236, 1276, 1630 and 1670). So this puts us
in \textbf{case 1(b)}.
\smallskip

\textbf{Case 1(b):}  There is a square $\alpha$ in $G_0$ and an
edge $e$ of $\alpha$ such that $G_0$ has a pair of
$\alpha$-opposite subgraphs with $e$ as the common edge of the
four squares linking $\g$ in both.

Assume, without loss of generality that $\g$ links 1234, and
123456 and 123478 are a pair of $\alpha$-opposite subgraphs that
$\g$ links with the same common edge say $\overline{12}$. So $\g$
links the square 1234, 1236, 1254, and 1256 in 123456,and 1234,
1238, 1274, and 1278 in 123478.  This forces $\g$ to link 1258,
because it links the three squares 1234, 1254, and 1238 in 123458.
 Similarly it forces $\g$ to link 1276 in the subgraph 123476.
Hence $\g$ links the squares 1234, 1254, 1236, 1256, 1238, 1258,
1274, 1276, and 1278. Thus $\g$ links every square in $G_0$
containing $\overline{12}$. The graph $G_0$ is isomorphic to
$K_{4,4}$, so each edge is contained in one component of a link of
two components \cite{Sac}. Take the non-split link $L_o\cup L_1$
in $G$ such that $L_o$ contains the edge $\overline{12}$.  Thus
$L_o$ links $\g$, and hence $(L\smallsetminus J)\cup L_o\cup L_1$
a non-split $n$-component link in $K_{2n+1,2n+1}$ and
$(L\smallsetminus J)\cup L_o$ is a non-split $(n-1)$-component
link.
\medskip

\textbf{Case 2:} The curve $\g$ links some $K_{3,3}$ subgraph of
$G$ in six squares.
\smallskip

Let $G_0$ be a subgraph of $G$ isomorphic to $K_{4,4}$ that
contains some $K_{3,3}$ subgraph that $\g$ links in six squares.
Label the two set of vertices of $G_0$ by \{1,3,5,7\} and
\{2,4,6,8\}.
\smallskip

\textbf{Case 2(a):}  There is no square $\alpha$ linking $\g$ in
$G_0$, such that the pair of $\alpha$-opposite subgraphs in $G$
both link $\g$ in six squares.
\smallskip

Without loss of generality $\gamma$ links $1234$ and links a
subgraph of $G_0$ isomorphic to $K_{3,3}$ containing 1234 in six
squares and the 1234-opposite subgraph in $G$ isomorphic to
$K_{3,3}$ in four squares.  Without loss of generality, $\gamma$
links 123456 in four squares, say, 1234, 1236, 1254, and 1256, and
links 123478 in six squares.  By Lemma \ref{lembig}, there are
three mutually nonadjacent edges in 123478 which each appear in
precisely two of six squares linking $\g$. Any square in 123478
will contain at least one of the edges in a set of three mutually
nonadjacent edges.  So each square in a set of six squares which
links $\g$ must contain exactly one of these three edges. Since
the set of six squares linking $\g$ in 123478 includes 1234,
precisely one of the mutually nonadjacent edges must be an edge of
1234. Thus the mutually nonadjacent edges are determined by taking
an edge of 1234 and two other nonadjacent edges of 123478 not in
the square 1234.  However, once the edge of 1234 is chosen, it
determines the remaining pair of mutually nonadjacent edges in
123478 which are not in 1234.  Since there are four choices for
the edge in 1234, there are four ways $\gamma$ can link 123478 in
six squares given that 1234 links $\gamma$. These four
possibilities are listed below with their three mutually
nonadjacent edges that appear exactly twice:
$$\begin{array}{cccc}
  \textrm{possibility 1} & \textrm{possibility 2} & \textrm{possibility 3} &
  \textrm{possibility 4} \\
  \overline{34}, \overline{27}, \overline{18} & \overline{23}, \overline{18}, \overline{47}
  & \overline{14}, \overline{27}, \overline{38} & \overline{12}, \overline{38}, \overline{47} \\
\end{array}$$
So the sets of squares that $\g$ links in 123478 are as follows:
$$\begin{array}{cccc}
  \textrm{possibility 1} & \textrm{possibility 2} & \textrm{possibility 3} &
  \textrm{possibility 4} \\
  1234 & 1234 & 1234 & 1234 \\
  1274 & 1274 & 3274 & 3274 \\
  1238 & 1438 & 1238 & 1438 \\
  1478 & 1278 & 3478 & 3278 \\
  3278 & 3278 & 1278 & 1278 \\
  3478 & 3478 & 1478 & 1478 \\
\end{array}$$
Possibilities 1, 2, and 3 contain 3478.  In these possibilities,
$\gamma$ links two disjoint squares in $G_0$, namely 3478 and
1256. Thus $(L\smallsetminus J)\cup 3478\cup 1258$ is a non-split
$n$-component link, and $(L\smallsetminus J)\cup 3478$ is a
non-split $(n-1)$-component link.  So we are done.

Now consider the possibility 4. At the beginning of this case we
assumed that $\g$ is linked to 1254, and since we are in
possibility 4, $\g$ also links 1278 and 1478.  There is no edge
that appears in all three of these squares so by Lemma
\ref{lembig}, $\g$ links 125478 in six squares.  Also, $\g$ links
123678 in 1236, 3278, and 1278. So, by Lemma \ref{lembig}, $\g$
also links 123678 in six squares. Thus $\g$ links both 125478 and
123678 in six squares including 1278.  Since 124578 and 123678 are
1278-opposite graphs, this violates the hypothesis of this case.
\smallskip

\textbf{Case 2(b):} There is a square $\alpha$ in $G_0$ linking
$\g$ such that the pair of $\alpha$-opposite subgraphs in $G_0$
both link $\g$ in six squares.
\smallskip

Without loss of generality, $\alpha =1234$, and the pair of
1234-opposite subgraphs 123456 and 123478 both link $\g$ in six
squares. Thus there are three mutually nonadjacent edges of
123456, each of which appears in two of the six squares that $\g$
links, and each of the other edges of 123456 appears in three of
the squares that $\g$ links. As we saw in \textbf{Case 2(a)},
precisely one of these three mutually nonadjacent edges appears in
each of the squares of 123456 linking $\g$. Since $\g$ links 1234,
without loss of generality $\overline{12}$ is one of these edges,
and hence $\overline{34}$ is not one of these edges.  Thus
$\overline{56}$ is also not one of these edges. Hence
$\overline{56}$ appears in three of the squares linking $\g$ in
123456.  Similarly, $\overline{78}$ appears in three of the
squares linking $\g$ in 123478.  In 123456 there are four squares
containing $\overline{56}$; they are of the form $ab56$ with
$a\in\{1,3\}$ and $b\in\{2,4\}$. There are also four squares
containing $\overline{78}$ in 123478 and they are of the form
$cd78$ with $c\in\{1,3\}$ and $d\in\{2,4\}$. Since three of each
of these sets of squares link $\g$ there must be two disjoint
squares $L_o$ and $L_1$ in $G_0$ that each link $\gamma$. Now
$(L\smallsetminus J)\cup L_o\cup L_1$ is a non-split $n$-component
link, and $(L\smallsetminus J)\cup L_o$ is a non-split
$(n-1)$-component link.
\end{proof}

\begin{rem}
For $K_{r,s}$ to be intrinsically $n$-linked it must contain $n$
disjoint simple closed curves.  The smallest simple closed curve
in a complete bipartite graph is a square.  A square contains four
vertices, two from each of the sets of vertices.  So the smallest
complete bipartite graph that could be intrinsically $n$-linked is
$K_{2n,2n}$.  Here we have shown that $K_{2n+1, 2n+1}$ is
intrinsically $n$-linked.  If $n=2$ it is known that $K_{4,4}$ is
intrinsically linked.  For $n>2$ there remain two graphs for which
it is not known whether $K_{2n,2n}$ or $K_{2n,2n+1}$ are
intrinsically $n$-linked.  Thus it is not known whether $K_{2n+1,
2n+1}$ is the smallest intrinsically $n$-linked bipartite graph.
\end{rem}

\begin{cor}
Let $n>1$, and $K_{2n+1,2n+1}$ be embedded in $\R^3$, every edge
of $K_{2n+1,2n+1}$ is contained in a non-split $n$-component link.
\end{cor}

\begin{proof}
Suppose there is an edge $e\in K_{2n+1,2n+1}$ that is not
contained in a non-split $n$-component link.  Then choose a
subgraph of $K_{2n+1,2n+1}$ isomorphic to $K_{4,4}$ that contains
the edge $e$. Since every edge of $K_{4,4}$ is contained in a
non-split 2-component link \cite{Sac}. The edge $e$ is contained
in a non-split 2-component link, $L_o\cup L_1$ in $K_{4,4}$.
Without loss of generality assume $e \in L_o$, and notice that
both $L_o$ and $L_1$ are squares. Now assume that $n>2$.  Next
choose a $K_{7,7}$ subgraph of $K_{2n+1,2n+1}$, which contains
this $K_{4,4}$. Then using the construction in the proof of the
Theorem \ref{Thm1}, with $J=L_1$ and $\g =L_o$ we obtain a
non-split 3-component link in $K_{7,7}$ containing $L_o$ and hence
the edge $e$. Thus we now assume that $n>3$.  The resulting link
is $L_o\cup L_2\cup L_3$, where either $L_2$ and $L_3$ link $L_o$
or $L_2$ links both $L_o$ and $L_3$. Now consider a $K_{9,9}$
subgraph of $K_{2n+1,2n+1}$, which contains $K_{7,7}$. In the
first case, where $L_2$ and $L_3$ link $L_o$, take $J=L_3$ and
$\g=L_o$ to get a non-split 4-component link in $K_{9,9}$
containing $L_o$. In the second case, where $L_2$ links $L_o$ and
$L_3$, take $J=L_3$ and $\g=L_2$ to get a non-split 4-component
link that contains $L_o$. In either case, the resulting link
contains $L_o$ and therefore it contains $e$. We can continue in
this way, at each stage choosing $J$ to be a square constructed in
the previous step and $\g$ to be a square linking $J$.  The
subsequent larger non-split link will always contain the edge $e$.
\end{proof}

With minor changes to the proof of Theorem \ref{Thm1} and
Corollary 1 we can prove a slightly stronger result about
tripartite graphs, as follows:

\begin{thm}\label{thm2}
Given $n>1$, every embedding of the complete tripartite graph
$K_{2n,2n,1}$ into $\mathbb{R}^3$ contains a non-split
$n$-component link.
\end{thm}

\begin{cor}
Let $K_{2n,2n,1}$ be embedded in $\R^3$, every edge of
$K_{2n,2n,1}$ is contained in a non-split $n$-component link.
\end{cor}

\section{Acknowledgment}
I would like to thank Erica Flapan for suggesting the subject of
intrinsically linked graphs and for many helpful conversations.

\end{document}